\documentclass[matzei,draft,numbook,envcountsame,envcountreset]{svjour}

\usepackage{graphics}

\newcommand{\er}{E_2(R)}
\newcommand{\ert}{E_2(R[t])}
\newcommand{\zz}{\vec{Z}}
\newcommand{\br}{B(R)}
\newcommand{\brt}{B(R[t])}
\newcommand{\fp}{\vec{F}_p}
\newcommand{\ez}{SL_2(\vec{Z})}
\newcommand{\ezt}{E_2(\vec{Z}[t])}
\newcommand{\slzt}{SL_2(\vec{Z}[t])}
\newcommand{\bz}{B(\vec{Z})}
\newcommand{\bzt}{B(\vec{Z}[t])}
\newcommand{\slfp}{SL_2(\vec{F}_p)}
\newcommand{\slfpt}{SL_2(\vec{F}_p[t])}
\newcommand{\bfp}{B(\vec{F}_p)}
\newcommand{\bfpt}{B(\vec{F}_p[t])}
\newcommand{\lra}{\longrightarrow}
\newcommand{\ra}{\rightarrow}
\newcommand{\bop}{\bigoplus}
\newcommand{\ot}{\otimes}
\newcommand{\op}{\oplus}

\begin{document}

\title{Amalgamated Free Products, Unstable Homotopy Invariance, and
the Homology of $SL_2(\vec{Z}[t])$}
\titlerunning{Amalgamated Free Products}
\author{Kevin P. Knudson\thanks{Supported by an NSF Postdoctoral Fellowship,
grant no. DMS-9627503}}

\institute{Department of Mathematics, Northwestern University, Evanston, IL
60208\\ \email{knudson@math.nwu.edu}}

\date{Received:  / Revised version:}

\maketitle

\begin{abstract}
We prove that if $R$ is a domain with many units, then the natural
inclusion $\er \ra \ert$ induces an isomorphism in integral homology.
This is a consequence of the existence of an amalgamated free product
decomposition of $\ert$.  We also use this decomposition to study the
homology of $\ezt$ and show that a great deal of the homology of
$\ezt$ maps nontrivially into the homology of $\slzt$. As a consequence,
we show that 
the latter is not finitely generated in all positive degrees.
\end{abstract}

\section*{Introduction}
\label{intro}

The fundamental theorem of algebraic $K$-theory asserts that if $A$ is
a regular ring, then there is a natural isomorphism $K_i(A[t])
\cong K_i(A)$ for all $i\ge 0$.  More generally, a
presheaf $\cal F$ on the category of schemes over a base $S$ 
is {\em homotopy invariant} if ${\mathcal F}(X\times
\vec{A}^1) \cong {\mathcal F}(X)$ for all schemes $X$.  The
utility of homotopy invariant presheaves is well-documented (see {\em
e.g.,} \cite{suslinicm,susvoe,voe}).

The existence of homotopy invariance in $K$-theory suggests that one
search for unstable analogues.  In a previous work \cite{knudson}, the
author showed that such a statement is true for infinite fields $k$:
the natural map $G(k)\ra G(k[t])$ induces an isomorphism in integral
homology for $G=SL_n,GL_n,PGL_n$.  Since homotopy invariance in
$K$-theory holds for all regular rings, one might hope that the above
isomorphism holds for rings other than infinite fields.

For a ring $A$, denote by $E_n(A)$ the subgroup of $GL_n(A)$ generated
by the elementary matrices.  In many cases, the group $E_n(A)$
coincides with $SL_n(A)$ ({\em e.g.,} $A$ local or Euclidean), but
this need not always be so.  In this paper, we study the homology of
the group $\ert$.  Our first result is the following.

\medskip

\noindent{\bf Theorem A}  {\em If $R$ is an integral domain with many
units, then the inclusion $\er\ra\ert$ induces an isomorphism
$$H_\bullet(\er,\zz)\stackrel{\cong}{\lra}H_\bullet(\ert,\zz).$$}
The definition of a ring with many units will be recalled
below.  The principal example of such a ring is a local ring with
infinite residue field (or an algebra over such a ring).

Let $k$ be an infinite
field and denote by $\tens{Sch}/k$ the category of smooth affine
schemes over $k$.  In this setting, Theorem A admits the following
geometric interpretation.

\medskip

\noindent{\bf Corollary} {\em The presheaves ${\mathcal H}_i =
H_i(E_2(-),\zz)$ are homotopy invariant on $\tens{Sch}/k$.}

\medskip

One would hope that a similar statement holds for $H_i(E_n(-),\zz)$
for $n>2$, but so far we have been unable to prove it.  What is clear,
as the following example of Weibel \cite{weibel} shows, is that if
Theorem A is to admit a generalization to the case $n>2$, then the
hypotheses on the ring $R$ will have to be strengthened.

\medskip

\noindent{\bf Example}  Consider the local ring of the node at the
singularity---$R=(k[x,y]/(y^2-x^3-x^2))_{(x,y)}$.  This is a
one-dimensional noetherian local domain, and, if $k$ is infinite, is a
ring with many units.  However, if $H_\bullet(E_n(R),\zz)\cong
H_\bullet(E_n(R[t]),\zz)$ for all $n$, then in particular
$H_2(E(R),\zz) \cong H_2(E(R[t]),\zz)$; {\em i.e.,} $K_2(R)\cong
K_2(R[t])$.  But in \cite{weibel} it is shown that $K_2(R)\ne
K_2(R[t])$.

\medskip

Thus, it seems that one must assume that $R$ is at least regular, and
it is not at all clear how to use this property.

There is another reason to study the homology presheaves ${\mathcal H}_i
=H_i(E_n(-),\zz)$.  If one could show that the ${\mathcal H}_i$ are 
homotopy invariant and admit transfer maps ({\em i.e.,} for a finite
flat morphism $R\ra S$ there is a map ${\mathcal H}_i(S) \ra
{\mathcal H}_i(R)$ satisfying the usual properties),
then one could prove the
Friedlander--Milnor conjecture \cite{fried} for $SL_n$:  If $G$
is a reductive group scheme over an algebraically closed field $k$ and
if $p\ne \textrm{char}\, k$, then $$H^\bullet_{et}
(BG_k,\zz/p)
\cong H^\bullet(BG(k),\zz/p).$$ 
Here, $BG_k$ is the simplicial classifying scheme of $G_k$ and
$BG(k)$ is the classifying space ($=$ simplicial set) of the
discrete group $G(k)$ of $k$-rational points. 
 Unfortunately, the ${\mathcal H}_i$ do not appear to
admit transfers.   However, it would still be 
interesting to know if the ${\mathcal H}_i$ are homotopy invariant.

Theorem A is a consequence of the existence of an amalgamated free
product decomposition $$\ert \cong \er *_{\br} \brt$$ where, for a
ring $A$, $B(A)$ denotes the upper triangular subgroup.  This follows
from Nagao's decomposition of $SL_2(k[t])$ ($k$ a field) and some
straight-forward group-theoretic considerations.  In particular, for
$R=\zz$ we have $\ezt \cong \ez *_{\bz} \bzt$ (note that $E_2(\zz) =
\ez$). This allows us to make partial computations of the homology of
$\slzt$.  Grunewald, Mennicke, and Vaserstein \cite{gmv} have shown
that $\slzt$ has free quotients $F$ of countable rank and that the
projection $\slzt \ra F$ may be taken so that all unipotent elements
(and hence all elements of $\ezt$) lie in the kernel.  This was
generalized and improved by Krsti\'c and McCool \cite{kmc}, who showed
that the same is true for any domain $R$ which is not a field.  This 
implies that $H_\bullet(\slzt,\zz)$ contains the
homology of a countably generated free group as a direct summand.
Note, however, that this only shows that $H_1(\slzt,\zz)$ has infinite
rank.

The amalgamated free product decomposition for $\ezt$ allows us to
compute its homology.  We show that $H_i(\ezt,\zz)$ contains a
countably generated free summand for $i\ge 1$.  It is not clear, a
priori, that any homology classes map nontrivially into
$H_\bullet(\slzt,\zz)$.  We show that this is indeed the case.

\medskip

\noindent{\bf Theorem B} {\em For each $i\ge 0$, the map
$$H_i(\ezt,\zz)\lra H_i(\slzt,\zz)$$ is nontrivial.}

\medskip

\noindent{\bf Corollary} {\em  For each $i\ge 1$, the group
$H_i(\slzt,\zz)$ is not finitely generated.}

\medskip

In the case $i=1$ we see that $H_1(\ezt)$ maps nontrivially into
$H_1(\slzt)$. Since all elementary matrices map to $1$ in
the free quotient $F$, the image of $H_1(\ezt)$ intersects the
copy of $H_1(F) \subset H_1(\slzt)$ trivially.  Hence we have
found a nontrivial, nonfinitely generated summand orthogonal
to $H_1(F)$.

Theorem B is proved by considering the (surjective!) homomorphisms
$\varphi_p:\ezt \ra \slfpt$ obtained by reducing polynomials modulo
$p$. We show that classes in $H_i(\ezt,\zz)$ map nontrivially into
$H_i(\slfpt,\zz)$ under $\varphi_{p*}$.  The result follows since
the map $\varphi_p$ factors through $\slzt$.

The paper is organized as follows.  In Section \ref{amal} we deduce
the amalgamated free product decomposition for $\ert$.  In Section
\ref{unstable} we prove Theorem A.  In Section \ref{ezt} we discuss
the homology of $\ezt$.  In Section \ref{fpt} we study the homology of
$\slfpt$.  Section \ref{compare} contains the proof of Theorem B.
Finally, Section \ref{other} deals with other homology classes in
$H_\bullet(\slzt,\zz)$.

\medskip

\noindent {\em Notation}  If $A$ is a ring (always assumed commutative
with unit), then $A^\times$ denotes the multiplicative group of units.
If $G$ is an abelian group, then ${\scriptstyle n}G$ denotes the
subgroup of $G$ annihilated by $n$.  If no coefficient module is
specified, homology groups are to be interpreted as integral homology.

\begin{acknowledgement}  I would like to thank Dick Hain and Andrei 
Suslin for many valuable conversations.  I am indebted to Chuck Weibel
for the example above.  Finally, I thank Mark Walker for listening
to me talk about homotopy invariance.
\end{acknowledgement} 

\section{Amalgamated Free Products}\label{amal}

Recall the following theorem of Nagao \cite{nagao}.  Let $k$ be a
field.  Then we have an amalgamated free product decomposition
$$SL_2(k[t]) \cong SL_2(k) *_{B(k)} B(k[t])$$ where for a ring $A$,
$B(A)$ denotes the subgroup of upper triangular matrices.  We will
make use of the following result of Serre (\cite{serre}, Prop.~3, 
p.~6).

\begin{proposition}\label{serreprop}  Suppose that $G$ is an
amalgamated free product, $G=G_1 *_A G_2$.  Let $H_i\subseteq G_i$ be
subgroups such that $H_1\cap A = B = H_2\cap A$.  Denote by $H$ the
subgroup of $G$ generated by the $H_i$.  Then the evident homomorphism
$$H_1 *_B H_2 \lra H$$ is an isomorphism. \hfill $\qed$
\end{proposition}

Let $R$ be an integral domain with field of fractions $Q$ and consider
the subgroups $\er \subseteq E_2(Q)$ ($=SL_2(Q)$) and $\brt \subseteq
B(Q[t])$.  Observe that $\er \cap B(Q) = \br = \brt \cap B(Q)$.

\begin{corollary}  There is a decomposition
$$\ert \cong \er *_{\br} \brt.$$
\end{corollary}

\begin{proof}  It is clear that the group $H$ generated by $\er$ and
$\brt$ is contained in $\ert$.  The reverse inclusion follows since
$\er$ contains the permutation matrix $\left(\begin{array}{cr}
                                           0  &  -1 \\
                                           1  &  0
                                       \end{array} \right)$
and hence $H$
contains all upper and lower triangular matrices. \hfill $\qed$
\end{proof}

\section{Rings with Many Units}\label{unstable}

In this section, we shall prove Theorem A.

\begin{definition}  A ring $A$ is an $S(n)$-ring if there exist
$a_1,\dots , a_n \in A^\times$ such that the sum of any nonempty
subfamily of the $\{a_i\}$ is a unit.  If $A$ is an $S(n)$-ring for every
$n$, then we say that $A$ has many units.
\end{definition}

Examples of rings with many units include local rings with infinite
residue fields and algebras over such rings.  However, a ring having
infinitely many units need not be a ring with many units.  An example
is the local ring $\zz_{(p)}$.  In fact, $\zz_{(p)}$ is an
$S(p-1)$-ring, but not an $S(p)$-ring.

Rings with many units satisfy the following property (see \cite{nessus},
Thm.~1.10).

\begin{theopargself}
\begin{proposition}\label{same}  If $R$ has
many units, then the inclusion $\br\ra \brt$ induces an isomorphism in
integral homology. \hfill $\qed$
\end{proposition}
\end{theopargself}

\begin{theorem}\label{thma}  Suppose $R$ is an integral domain with
many units.  Then the inclusion $\er\ra\ert$ induces an isomorphism
$$H_\bullet(\er,\zz)\stackrel{\cong}{\lra} H_\bullet(\ert,\zz).$$
\end{theorem}

\begin{proof}  The amalgamated free product decomposition yields a
Mayer--Vietoris sequence for computing $H_\bullet(\ert,\zz)$.  Since
the evident map $H_\bullet(\br,\zz)\ra H_\bullet(\brt,\zz)$ is split injective
for any domain $R$, the long exact sequence breaks up into short exact
sequences
$$0\ra H_i(\br) \ra H_i(\brt)\op H_i(\er) \ra H_i(\ert) \ra 0$$
(this is valid for any domain $R$).  By Proposition \ref{same},
$H_i(\br)\cong H_i(\brt)$ if $R$ has many units so that the map
$$H_i(\er,\zz)\lra H_i(\ert,\zz)$$ is an isomorphism. \hfill $\qed$
\end{proof}

\begin{remark}  Evidently, this approach will not work for computing
$H_\bullet(E_n(R[t]))$ for $n\ge 3$ since no such amalgamated free
product decomposition exists.
\end{remark}

\section{The Homology of $\ezt$}\label{ezt}

We now shift gears and study the homology of $\ezt$.  Theorem
\ref{thma} does not apply in this case since $\zz$ does not have many
units.  It is possible to calculate
$H_\bullet(\ezt,\zz)$ explicitly, but the final answer is a bit
complicated.  For this reason, we will compute only a part of
the integral homology in detail.  We also compute the homology of
$\ezt$ with $\fp$-coefficients.

For each $i$ we have a short exact sequence
$$0\ra H_i(\bz)\ra H_i(\bzt)\op H_i(\ez) \ra H_i(\ezt) \ra 0.$$
Observe that $\bzt = \bz \times t\zz[t]$ so that for any principal
ideal domain $k$ we have the K\"unneth exact sequence
\begin{eqnarray*}
\lefteqn{0\ra \bop_{l+m=i}H_l(\bz,k)\ot H_m(t\zz[t],k) \ra H_i(\bzt,k) \ra} 
\hspace{1.25in} \\
& & \bop_{l+m=i-1} 
\textrm{Tor}^k_1(H_l(\bz,k),H_m(t\zz[t],k))\ra 0.
\end{eqnarray*}
In particular, if $k$ is a field we have
$$H_\bullet(\bzt,k) \cong H_\bullet(\bz,k)\ot H_\bullet(t\zz[t],k).$$

Suppose $k$ is the finite field $\fp$.  Then we have the following
result.

\begin{proposition}\label{bhomfp} If $p\ge 3$, then for
$R=\zz,\zz[t]$, $$H_\bullet(\br,\fp) = H_\bullet(R,\fp) = 
{\bigwedge\nolimits}_{\fp}^\bullet (R\ot \fp).$$  If $p=2$, then
$$H_\bullet(\bzt,\vec{F}_2) = (H_\bullet(\zz/2,\vec{F}_2)\ot
H_\bullet(\zz,\vec{F}_2)) \ot {\bigwedge\nolimits}_{\vec{F}_2}^\bullet
t\vec{F}_2[t].$$
\end{proposition}

\begin{proof}  Recall that if $A$ is an abelian group then there is an
isomorphism $${\bigwedge\nolimits}_{\fp}^\bullet (A\ot\fp) \ot
\Gamma_{\fp}({\scriptstyle p}A) \stackrel{\cong}{\lra}
H_\bullet(A,\fp)$$ where $\Gamma_{\fp}$ is a divided power algebra
\cite{brown}, p.~126.  The isomorphism is natural if $p\ne 2$.  The
first result follows from the fact that $\br = \zz/2 \times R$
($R=\zz,\zz[t]$) and that $H_i(\zz/2,\fp)$ vanishes for $i>0$ if $p\ge
3$.  Note that in this case ${\scriptstyle p}R = 0$.  The second
assertion is just a restatement of the K\"unneth formula. \hfill $\qed$
\end{proof}

\begin{corollary}\label{e2pge3}  If $p\ge 3$, the homology of $\ezt$
is $$H_i(\ezt,\fp) = \left\{ \begin{array}{ll}
                         \fp & i=0 \\
                        t\fp[t] \op (\zz/12 \ot\fp) & i=1 \\
                     \bigwedge_{\fp}^i \fp[t] \op
                         H_i(\ez,\fp) & i\ge 2.
                     \end{array} \right. $$
In particular, if $p\ge 5$, then for $i\ge 2$
$$H_i(\ezt,\fp) = {\bigwedge\nolimits}_{\fp}^i \fp[t].$$
\end{corollary}

\begin{proof}  Consider the short exact sequence
\begin{eqnarray*}
\lefteqn{0\ra H_i(\bz,\fp)\ra H_i(\bzt,\fp)\op H_i(\ez,\fp) \ra} 
\hspace{2.75in} \\
&  & H_i(\ezt,\fp) \ra 0.
\end{eqnarray*}
The result for $i\ge 2$ follows easily since 
$$H_i(\bz,\fp) = H_i(\zz,\fp) = 0.$$  For $i=1$, we have
$H_1(\bz,\fp) = \fp$ and $H_1(\bzt,\fp) = \fp[t].$  The map $\fp \ra
\fp[t]$ is the usual inclusion and the calculation of $H_1(\ezt,\fp)$
follows easily.  The final assertion is a consequence of the fact that
the abelianization map $\ez \ra \zz/12$ induces an isomorphism on
integral homology. \hfill $\qed$
\end{proof}

The calculation for $p=2$ is more complicated, and instead of writing
down the complete answer we make the following observation.

\begin{proposition}\label{e2p2}  For each $i$, $H_i(\ezt,\vec{F}_2)$
contains $\bigwedge_{\vec{F}_2}^i t\vec{F}_2[t]$ as a
direct summand.
\end{proposition}

\begin{proof}  Since $\bzt = \bz \times t\zz[t]$, it follows that
the quotient $H_i(\bzt,\vec{F}_2)/H_i(\bz,\vec{F}_2)$ consists of
$$\bop_{p+q=i, p<i} H_p(\bz,\vec{F}_2) \ot H_q(t\zz[t],\vec{F}_2),$$ 
and this contains $H_0(\bz,\vec{F}_2) \ot
H_i(t\zz[t],\vec{F}_2)= \bigwedge_{\vec{F}_2}^i
t\vec{F}_2[t]$. \hfill $\qed$
\end{proof}

The integral homology is similarly complicated and we leave it to the
interested reader to write it down. What is important for our purposes is
the following.

\begin{proposition}\label{e2int}  For each $i$, $H_i(\ezt,\zz)$ contains
$\bigwedge_{\zz}^i t\zz[t]$ as a direct summand.  Moreover, the element
$t^{l_1} \wedge \cdots \wedge t^{l_i}$ maps to the element 
$t^{l_1} \wedge \cdots \wedge t^{l_i}$ in $\bigwedge_{\fp}^i t\fp[t]
\subseteq H_i(\ezt,\fp)$ under the map induced by reducing coefficients
modulo $p$.
\end{proposition}

\begin{proof}  Arguing as in the proof of Proposition \ref{e2p2} we see
that the group $H_i(\bzt,\zz)/H_i(\bz,\zz)$ contains a copy of the
tensor product $H_0(\bz,\zz) \ot H_i(t\zz[t],\zz)$.  Since $t\zz[t]$
is a torsion-free abelian group, its integral homology is simply an
exterior algebra; {\em i.e.,} $H_i(t\zz[t],\zz) = \bigwedge_{\zz}^i
t\zz[t]$.
The second assertion follows by considering the homomorphism
of Mayer--Vietoris sequences induced by reducing coefficients modulo
$p$. \hfill $\qed$
\end{proof}

\section{The Homology of $\slfpt$}\label{fpt}

In this section we compute the homology of the group $\slfpt$.  As
before, the amalgamated free product decomposition of $\slfpt$ yields
a collection of short exact sequences
\begin{eqnarray*}
\lefteqn{0 \ra H_i(\bfp) \ra H_i(\bfpt)\op H_i(\slfp) \ra} \hspace{2.75in} \\
&  & H_i(\slfpt) \ra 0.
\end{eqnarray*}
  Observe that if $\ell$ is a prime distinct from $p$, then the
natural map $$H_i(\bfp,\vec{F}_\ell) \lra H_i(\bfpt,\vec{F}_\ell)$$ is
an isomorphism for all $i$.  This follows by considering the Hochschild-
Serre spectral sequence associated to the extension
\begin{equation}\label{ext}
0 \lra R \lra \br \lra \fp^\times \lra 1
\end{equation}
for $R=\fp,\fp[t]$.  Since $H_q(R,\vec{F}_\ell) = 0$ for $q>0$, the 
claimed isomorphism follows.  As a consequence, we see that for all $i$,
$$H_i(\slfpt,\vec{F}_\ell) = H_i(\slfp,\vec{F}_\ell)$$ for $\ell \ne p$.
(This was proved, in greater generality, by C. Soul\'e \cite{soule}.)

If we consider integral or $\fp$ coefficients, however, the situation is
much different.

\begin{lemma}\label{action}  For $R=\fp,\fp[t]$, we have
$$H_i(\br,\fp) = H_0(\fp^\times,H_i(R,\fp)).$$
\end{lemma}

\begin{proof}  Consider the Hochschild--Serre spectral sequence
associated to the extension (\ref{ext}); it has $E^2$-term
$$E^2_{r,s} = H_r(\fp^\times, H_s(R,\fp)).$$  Since $|\fp^\times|
=p-1$ is invertible in $\fp$, we see that $E^2_{r,s}=0$ for $r>0$
(see \cite{brown}, p.~84). \hfill $\qed$
\end{proof}

The action of $\alpha \in \fp^\times$ on $R$ is $\alpha:x \mapsto
\alpha^2x$.  Thus, if $p=2\; \textrm{or}\; 3$, we have the following.

\begin{corollary}  If $p=2$ or $3$, then for all $i\ge 0$,
$$H_i(\br,\fp) = H_i(R,\fp).$$
\end{corollary}

\begin{proof}  In this case, $(\vec{F}_2^\times)^2 = \{1\} =
(\vec{F}_3^\times)^2$ so that $\fp^\times$ acts trivially on
$H_i(R,\fp)$. \hfill $\qed$
\end{proof}

\begin{corollary}\label{sl2p23} If $p=2$ or $3$, then for all $i\ge 0$,
\begin{eqnarray*}
H_i(\slfpt,\fp) & = & H_i(\slfp,\fp) \\
 & & \op H_i(\bfpt,\fp)/H_i(\bfp,\fp).
\end{eqnarray*}
In particular, $H_i(\slfpt,\fp)$ contains a copy of
$\bigwedge_{\fp}^i t\fp[t]$ as a direct summand.
\end{corollary}

\begin{proof}  If $p=2$ or $3$, then $\bfpt = \bfp \times t\fp[t]$.  Hence
the group $H_i(\bfpt,\fp)/H_i(\bfp,\fp)$ consists of
$$\bop_{l+m=i, l<i} H_l(\fp,\fp)\ot H_m(t\fp[t],\fp).$$ 
This contains $H_i(t\fp[t],\fp)$ as a summand and this, in turn, 
contains $\bigwedge_{\fp}^i t\fp[t]$ as a summand. \hfill $\qed$
\end{proof}

The calculation of the integral homology is obviously more complicated,
even for $p=2,3$.  However, we make the following observation.

\begin{proposition}\label{wedgefp}  If $p=2$ or $3$, then for all $i\ge 0$,
$H_i(\slfpt,\zz)$ contains a copy of $\bigwedge_{\zz}^i t\fp[t]$. 
Furthermore, the element $t^{l_1}\wedge \cdots \wedge t^{l_i}
\in \bigwedge_{\zz}^i t\fp[t]$ maps to $t^{l_1}\wedge \cdots \wedge t^{l_i}
\in \bigwedge_{\fp}^i t\fp[t] \subseteq H_i(\slfpt,\fp)$ under the map
induced by reducing coefficients modulo $p$.
\end{proposition}

\begin{proof}  For any abelian group $A$, the map $\bigwedge_{\zz}^
\bullet A \ra H_\bullet(A,\zz)$ is injective \cite{brown}, p.~123.
It follows that $H_i(\bfpt,\zz)$ contains a copy of $H_0(\bfp,\zz)
\ot H_i(t\fp[t],\zz)$, which, in turn, contains a copy of
$\bigwedge_{\zz}^it\fp[t]$.
Thus, the group $H_i(\slfpt,\zz)$ does also.  The second assertion
follows from the naturality of the map induced in homology by the
coefficient homomorphism $\zz\ra\fp$. \hfill $\qed$
\end{proof}

The situation for primes greater than $3$ is made more difficult by
the fact that $(\fp^\times)^2 \ne \{1\}$ for $p\ge 5$.  Thus, the
calculation of $H_\bullet(\bfp)$ and $H_\bullet(\bfpt)$ is more
involved.  Since considering the primes $2$ and $3$ is sufficient
for our purposes here, we leave the case $p\ge 5$ to the interested
reader.

\section{The Homology of $\slzt$}\label{compare}

In this section we prove Theorem B (which was stated in the introduction).
We first demonstrate the following result.

\begin{proposition}\label{reduce}  If $p=2$ or $3$, then the natural map
$$H_i(\ezt,\zz)\ra H_i(\slfpt,\zz)$$ induced by $\varphi_p:\ezt\ra\slfpt$
maps the element $t^{l_1}\wedge \cdots \wedge t^{l_i} \in
\bigwedge_{\zz}^i t\zz[t]$ to  $t^{l_1}\wedge \cdots \wedge t^{l_i} \in
\bigwedge_{\zz}^i t\fp[t]$.  The same is true with $\fp$-coefficients.
\end{proposition}

\begin{proof}  Consider the commutative diagram
$$\begin{array}{ccccc}
H_i(\bz) & \ra & H_i(\bzt)\op H_i(\ez) & \ra &
H_i(\ezt) \\
 \downarrow &    &     \downarrow        &     &
\downarrow \\
H_i(\bfp) &\ra & H_i(\bfpt)\op H_i(\slfp) & \ra &
H_i(\slfpt). 
\end{array}$$  Since $t^j \in t\zz[t]$ maps to $t^j \in t\fp[t]$,
the claim is clear.  The second assertion follows similarly. \hfill
$\qed$
\end{proof}

Consider the inclusion $j:\ezt\ra \slzt$ and the induced map $j_*$
on homology.  For each $i$-tuple $\underline{\ell} = (l_1,\dots ,l_i)$,
$l_1<\cdots <l_i$, denote by $x_{\underline{\ell}}$ the homology
class $j_*(t^{l_1}\wedge \cdots \wedge t^{l_i}) \in
H_i(\slzt,\zz)$.

\begin{theorem}\label{thmb}  For each $i$-tuple $\underline{\ell}$, the class
$x_{\underline{\ell}}$ is nontrivial.  
Moreover, if $\underline{m}$ is another $i$-tuple
distinct from $\underline{\ell}$, then $x_{\underline{\ell}}$ and 
$x_{\underline{m}}$ are
distinct.  Thus, the group $H_i(\slzt,\zz)$ is not finitely generated
for all $i\ge 1$.
\end{theorem}

\begin{proof} Let $p=2$ or $3$. Consider the commutative diagram
$$\begin{array}{ccc}
\ezt &  \stackrel{j}{\lra}       & \slzt  \\
     & \varphi_p \searrow           & \downarrow \pi_p \\
     &                           &\slfpt
\end{array}$$
and the induced commutative diagram
$$\begin{array}{ccc}
H_i(\ezt,\zz)  &  \stackrel{j_*}{\lra}   &  H_i(\slzt,\zz) \\
               & \varphi_{p*} \searrow      &  \downarrow \pi_{p*} \\
               &                         & H_i(\slfpt,\zz)
\end{array}$$
in homology.  Since $\varphi_{p*}(t^{l_1}\wedge \cdots \wedge t^{l_i})
= t^{l_1}\wedge\cdots\wedge t^{l_i}$, we see that 
$\pi_{p*}j_*(t^{l_1}\wedge\cdots\wedge t^{l_i}) \ne 0$ and hence
$x_{\underline{\ell}}$ is nontrivial.  Moreover, since $\varphi_{p*}(t^{l_1}
\wedge\cdots\wedge t^{l_i}) \ne \varphi_{p*}(t^{m_1}
\wedge\cdots\wedge t^{m_i})$, we see that $x_{\underline{\ell}}\ne 
x_{\underline{m}}$ provided $\underline{\ell} \ne \underline{m}$.
\hfill $\qed$
\end{proof}

\begin{corollary}\label{order6}  For each $i$-tuple $\underline{\ell}$, the
element $x_{\underline{\ell}}$ has either infinite order or order divisible by
$6$.
\end{corollary}

\begin{proof} Note that $\varphi_{p*}(t^{l_1}\wedge\cdots\wedge t^{l_i})$ is
an element of order $p$ in the group $H_i(\slfpt,\zz)$ (since this group is all
$p$-torsion).  It follows that $x_{\underline{\ell}}$ has order divisible by
the least common multiple of $2$ and $3$ (if it is not infinite). \hfill
$\qed$
\end{proof}

\begin{remark} We can glean further information by reducing modulo
primes $p\ge 5$.  By Lemma \ref{action}, $$H_i(\bfpt,\fp) = H_0(\fp^\times,
H_i(\fp[t],\fp)),$$ where $\alpha \in \fp^\times$ acts on $\fp[t]$
via $\alpha:x\mapsto \alpha^2x$.  Now, $H_i(\fp[t],\fp)$ contains
$\bigwedge_{\fp}^it\fp[t]$ as a direct summand and if $i$ is a
multiple of $(p-1)/2$, then, by Fermat's Little Theorem, $\fp^\times$
acts trivially on $\bigwedge_{\fp}^it\fp[t]$.  This implies that for
$i=n(p-1)/2$, $H_i(\slfpt,\fp)$ contains a copy of $\bigwedge_{\fp}^i
t\fp[t]$.  Hence, if $\underline{\ell}$ is an $i$-tuple, $i=n(p-1)/2$,
then $x_{\underline{\ell}}$ has order divisible by $2\cdot 3\cdot p$.
For example, the $x_{\underline{\ell}}$ in $H_{2n}$ have order divisible
by $2\cdot 3\cdot 5=30$, those in $H_{3n}$ have order divisible by 
$2\cdot 3\cdot 7=42$, etc.  Given this,
it seems reasonable to conjecture that the $x_{\underline{\ell}}$ have
infinite order. 
\end{remark}

\section{Other Classes in $H_\bullet(\slzt,\zz)$}\label{other}

In \cite{kmc}, Krsti\'c and McCool exhibit a basis of a free retract
$F$ of the quotient $\slzt/U_2(\zz[t])$ (here, $U_2(\zz[t])$ denotes the
subgroup generated
by unipotent matrices).  It consists of elements $h_{p,k}$ where $p\ge 2$
and $k\ge 1$.  The elements $h_{p,k}$ are defined as
$$h_{p,k} = \left(\begin{array}{cc}
                   1 + pt^k  &  t^{3k}    \\
                      p^3    & 1-pt^k + p^2t^{2k}
                   \end{array} \right).$$
Thus, the elements $\overline{h}_{p,k} \in H_1(\slzt,\zz)$ form the 
basis of a free direct summand; denote this summand by 
$\zz\{\overline{h}_{p,k}\}$.  Then we may write
$$H_1(\slzt,\zz) = H_1(\ez,\zz) \op \zz\{\overline{h}_{p,k}\} 
\op X,$$ where $X$ is some abelian group.

Since $\ezt \subset U_2(\zz[t])$, the elements 
$$x_k = \left(\begin{array}{cc}
                 1   &  t^k  \\  0  &  1
               \end{array} \right)$$
map to $1$ in $F$.  Thus, the elements $\overline{x}_k$ map into
the summand $X \subset H_1(\slzt,\zz)$.  Hence we see that $X$ is
also not finitely generated.  

Consider the homomorphism $\pi_p:\slzt\ra\slfpt$.  Then we have the
following.

\begin{proposition} The kernel of the map $$\pi_{p*}:H_1(\slzt,\zz)\ra
H_1(\slfpt,\zz)$$ is not finitely generated.
\end{proposition}

\begin{proof}  If $p\ge 5$, then $H_1(\slfpt,\zz) = H_1(\slfp,\zz)=0$ so
that $\pi_{p*}=0$.
If $p=2$ or $3$ then it is not clear that any summand lies in the
kernel.  But note that $\pi_p(x_k) = \pi_p(h_{p,k})$ so that
$\pi_{p*}(\overline{x}_k - \overline{h}_{p,k}) = 0$ for all $k\ge 1$.
\hfill $\qed$
\end{proof}

Consider the elements $g_{p,k}$ defined by
$$g_{p,k} = \left(\begin{array}{cc}
                     1   &   -t^k   \\
                     -p      &  1+pt^k
                  \end{array}\right).$$
Note that each $g_{p,k}$ is not unipotent; {\em i.e.,} the matrix
$$n_{p,k} = \left(\begin{array}{cc}
                     0   &   -t^k    \\
                     -p    &   pt^k
                  \end{array}\right)$$
has infinite order.  However, we do have the following.

\begin{lemma}  For each $k,l$, $g_{p,k} \equiv g_{p,l}$ modulo
$U_2(\zz[t])$.
\end{lemma}

\begin{proof}  An easy calculation shows that
$$g_{p,k}^{-1}g_{p,l} = \left(\begin{array}{cc}
                                  1  &  t^k - t^l  \\
                                  0  &     1
                              \end{array}\right).$$
The result follows. \hfill $\qed$
\end{proof}

Thus, the set of all $g_{p,k}$ (for a fixed $p$) lies in a single
coset of $U_2(\zz[t])$.  Denote by $\overline{g}_p$ the induced (nontrivial!)
element of $\zz\{\overline{h}_{p,k}\}$.  

\begin{proposition} If $p=2$ or $3$, then $\overline{g}_p +
\overline{x}_k$ lies in the kernel of $\pi_{p*}$ for all $k\ge 1$.
\end{proposition}

\begin{proof} Note that $\pi_p(g_{p,k}) = x_k^{-1}$.  It follows
then that $\pi_{p*}(\overline{g}_p) = -t^k \in \bigwedge_{\zz}^1
t\fp[t] \subset H_1(\slfpt,\zz)$.  Since $\pi_{p*}(\overline{x}_k)
= t^k$, the result follows. \\ \mbox{} \hfill $\qed$
\end{proof}

Note that we also have $\overline{g}_p + \overline{h}_{p,k} 
\in \textrm{ker}(\pi_{p*})$ since $\overline{g}_p + \overline{x}_k$
and $\overline{h}_{p,k} - \overline{x}_k$ lie in the kernel.

\end{document}